\magnification=1200

{\bf Holomorphic factorization of matrices of polynomials}

\bigskip
John P. D'Angelo

Dept. of Mathematics

University of Illinois  

Urbana IL 61801 

USA

\bigskip

{\bf Introduction}
\bigskip

This paper considers some work done
by the author and Catlin [CD1,CD2,CD3] concerning positivity conditions
for bihomogeneous polynomials and metrics on bundles
over certain complex manifolds.
It presents a simpler proof of a special case of the main result in [CD3], 
providing also
a self-contained proof of a generalization of the main result from [CD1].
Some new examples and applications appear here as well. 
The idea is to use the Bergman
kernel function and some operator theory to prove purely algebraic
theorems about matrices of
polynomials.

The main idea arises from generalizing
a classical factorization question. See [Dj] and [RR] for many
aspects of factorization of non-negative matrices and operators on
Hilbert spaces.
Consider a real-analytic
matrix-valued function $F(z,{\overline z})$
that is positive semi-definite at each point. Is there a holomorphic
matrix-valued function $A(z)$ such that $F(z,{\overline z})= A(z)^* A(z)$?
Here $A^*$ denotes the conjugate transpose of $A$.
In general the answer is no, even when $F$ is a scalar, positive-definite,
and its entries are
bihomogeneous polynomials. Because such factorizations 
have many applications,
we allow ourselves a generalization; 
we can first multiply $F$ by powers of a
scalar function $R$, and ask whether we can factor $R^dF$ for 
sufficiently large $d$. 
This is a natural thing to
do when one studies proper holomorphic mappings between balls in 
different dimensions, and one chooses $R$ to be the squared Euclidean norm.
See [CD1] for applications. This multiplication
also admits an interpretation in terms
of metrics on tensor products of Hermitian line bundles.

We write $\langle \zeta, w \rangle $ for the Euclidean Hermitian inner
product of $\zeta$ and $w$ on any
finite-dimensional complex Euclidean space, and $||\zeta||^2$ for the
squared norm. Later we will use subscripts to denote
$L^2$ norms. We write ${\bf V}(A)$ to denote the variety defined by the
simultaneous vanishing of the component functions of a holomorphic
mapping $A$.

A {\it bihomogeneous polynomial} on ${\bf C^n}$ is a polynomial function
$f : {\bf C^n} \times {\bf C^n} \to {\bf C}$ that is homogeneous
of the same degree $m$ in each set of variables. 
We will be considering $f(z,{\overline w})$, which is conjugate-analytic in
the second set of variables. 
The polynomial defined by $f(z,{\overline z})$ is
real-valued if and only if 
$f(z,{\overline w}) = {\overline {f(w,{\overline z})}}$;
we call such an 
$f$ a bihomogeneous real-valued polynomial
on ${\bf C^n}$ of degree $2m$.

Suppose that $f$ is a bihomogeneous real-valued polynomial of degree $2m$.
By elementary linear algebra, it is possible
to find holomorphic polynomial mappings $A$ and
$B$, with finitely many components, that are homogeneous of degree $m$, and 
such that

$$ f(z,{\overline w}) = \langle A(z),A(w) \rangle - \langle B(z),B(w) \rangle \eqno (1)$$
It follows from (1) that
$$ f(z,{\overline z}) = || A(z)||^2 - ||B(z)||^2 \eqno (2) $$
Suppose in addition
that $f(z,{\overline z}) \ge 0$. We investigate the following questions.
Can we choose $B=0$ in (2), and if we cannot, can we do this for
 $R^d f$, where $R$ is an appropriate multiplier and $d$ is sufficiently large? 
Suppose more generally that
$F(z,{\overline z})$ is a matrix of bihomogeneous polynomials each of degree $2m$,
and that is positive semi-definite at each $z \ne 0$.
Can we factor $||z||^{2d} F(z,{\overline z})$ for sufficiently large $d$?

The following theorem gives a decisive answer 
in the positive-definite case.

\bigskip

{\bf Theorem 1}. [Catlin-D'Angelo]. Suppose that $f$ is 
a bihomogeneous real-valued polynomial on ${\bf C^n}$ of degree $2m$.
 Then $f$ is positive
away from the origin if
and only there is an integer $d$ and a holomorphic homogeneous polynomial
mapping $A$, whose components span the space of homogeneous polynomials of
degree $m+d$, such that
$$ ||z||^{2d} f(z,{\overline z}) = ||A(z)||^2.  \eqno (3) $$
Suppose that $F(z,{\overline z})$ is an $r$ by $r$ matrix whose
entries are bihomogeneous polynomials of degree $2m$. 
Then $F(z,{\overline z})$ 
is positive-definite at each point $z \ne 0$ if and only if
there is an integer $d$ and a holomorphic homogeneous polynomial
matrix $A$, whose row vectors span the space of $r$-tuples of homogeneous
polynomials of degree $m+d$, such that
$$ ||z||^{2d} F(z,{\overline z}) = A(z)^* A(z). \eqno (4) $$

Note that (3) is the scalar
version of (4).
The scalar statement about $f$ was proved in [CD1]. The matrix version
is a special case of
a general result from [CD3] about Hermitian metrics on bundles over certain
complex manifolds. Because it is a special case, some steps in the proof simplify;
its intrinsic interest justifies giving the simpler proof here. 
In the scalar statement, one can
replace the condition that the components
of $A$ span, by the condition that ${\bf V}(A) = \{0\}$; the exponent required 
may be smaller. The proof reveals that the stronger condition on
$A$ is the natural one.

When $F(z,{\overline z}) = A(z)^* A(z)$, necessarily
$F(z,{\overline z})$ is positive semi-definite at each point.
A general result such as Theorem 1 cannot
hold in the positive semi-definite case, for the following simple reason.
Suppose $r=1$ for simplicity.
If the zero set of $f$ is not an analytic variety, then 
there is no hope to write
$ ||z||^{2d} f(z,{\overline z}) = ||A(z)||^2 $, because 
the zero set of the right side
is an analytic variety. A simple example where the zero set
of a bihomogeneous polynomial
fails to be an analytic variety is given by (5).
$$ f(z,{\overline z}) = (|z_1|^2 - |z_2|^2)^2 . \eqno (5) $$
Some results hold under specific hypotheses in
the semi-definite case, but most of this paper considers only the
positive-definite case.

We give a complete proof of Theorem 1, relying on the Bergman kernel
function for the unit ball. In Theorem 2 we give a simple application
to elliptic PDE. In Theorem 3 we 
reinterpret Theorem 1 in terms of the universal bundle over
complex projective space. We also provide some illuminating examples along the way. 
We close the paper with some brief remarks about factorization theorems proved in
the 1970s.

\bigbreak

{\bf I. Holomorphic factorability}
\bigskip

Suppose that we are given an $r$ by $r$ matrix $F(z,{\overline w})$ whose
entries $F_{ij}(z,{\overline w})$ are bihomogeneous polynomials of degree $2m$
on ${\bf C^n}$.
Let $N=N(n,m)$ denote the dimension of the vector space  
of homogeneous polynomials of degree $m$
in $n$ variables. We write $V_d$ for the vector space of $r$-tuples of homogeneous
polynomials of degree $d$ on ${\bf C^n}$.

{\bf Definition 1}. The $r$ by $r$
matrix of bihomogeneous polynomials $F(z,{\overline w})$ 
of degree $2m$ is called 
{\it holomorphically factorable} if there is an integer $s$ 
and a matrix $(E_{jk}(z))$, for $ j=1,...,r$ and $k=1,...,s$, 
of homogeneous
polynomials of degree $m$ such that
(6) holds.

$$ F_{ij} (z,{\overline w}) = \langle E_j(z), E_i(w) \rangle = 
\sum_{k=1}^s E_{jk}(z){\overline {E_{ik}(w)}}  \eqno (6) $$
Here $E_j(z) = (E_{j1}(z),...,E_{js}(z)) $.

Let $A$ be the transpose of $E$;
we observe immediately that Definition 1 implies that
$$ F(z,{\overline w}) = A(w)^* A(z). \eqno (7) $$

{\bf Definition 2}. The matrix of bihomogeneous polynomials $F(z,{\overline w})$ 
of degree $2m$ is called 
{\it strictly holomorphically factorable} if it is holomorphically factorable, and in
addition the $s$ column vectors of $(E_{jk}(z))$ are a basis for 
$V_d$.

Note that the notion of strict holomorphic factorability requires that we 
choose $s=N$. The first concept allows us to write (7), forcing
$F(z,{\overline z}) $ to be positive semi-definite. The second concept 
ensures that $F(z,{\overline z})$ is positive-definite, 
but it is an even stronger condition. 
We give a simple example when $r=1$
to illustrate the difference.

{\bf Example}. Let $F(z,{\overline w})$ be the one-by-one matrix given by
$z_1^2 {\overline w_1}^2 + z_2^2 {\overline w_2}^2 $. 
Here Definition 1 holds with
$E_1(z) = (z_1^2, z_2^2)$.
Note that $F(z,{\overline z})$ is positive-definite away
from the origin, but the absence of the $z_1z_2$ 
cross-term means that Definition 2 fails.

\bigskip

{\bf II. The link to operator theory}
\bigskip

We will prove Theorem 1 by using some facts about compact operators on
a Hilbert space. Because we consider (matrices of)
bihomogeneous polynomials and use the
Euclidean norm as a multiplier, it suffices to consider
the Hilbert space ${\cal H}$ of $r$-tuples of $L^2$ functions on the unit ball
$B$ in
complex Euclidean space ${\bf C^n}$. We let ${\cal A}$ denote the closed subspace of
$r$-tuples of holomorphic functions in $L^2(B)$. We write $P$ for the Bergman
projection from ${\cal H}$ to ${\cal A}$; it is the usual Bergman projection acting
on each component.

There is an orthogonal decomposition of ${\cal A}$ into 
the finite-dimensional subspaces
$V_d$. 
The following result
provides the crucial link between matrices whose entries are bihomogeneous
polynomials and operator
theory on ${\cal A}$.

{\bf Proposition}. Let $F(z,{\overline z})= (F_{ij}(z,{\overline z}))$ 
be an $r$ by $r$
 matrix of bihomogeneous
polynomials of degree $2d$ on ${\bf C^n}$. The following are equivalent:

1. The matrix $F$ is {\it strictly holomorphically factorable}.
Thus there is an $s$ by $r$  matrix $E(z)$ of holomorphic homogeneous polynomials 
whose column vectors give a basis for $V_d$ 
for each $z \ne 0$. Furthermore, with $A$ the transpose of $E$,  (8) holds.

$$ F(z,{\overline w}) = A(w)^* A(z)  \eqno (8) $$

2. Consider the integral operator $ T: V_d \to V_d $ defined by

$$ Th(z)_j = \int_B \sum_i F_{ij} (z,{\overline w}) h_i(w) dV(w)  \eqno (9) $$
Then $T$ is positive on $V_d$. That is,
there is a positive constant $c$ so that for $h \in V_d$, 

$$ \langle Th,h \rangle_{\cal H} \ge c ||h||^2_{\cal H}  \eqno (10) $$

3.  Write 
$$ F_{ij}(z,{\overline w}) = \sum F_{ij\alpha\beta} z^\alpha {\overline w}^\beta
\eqno (11) $$
for constants $ F_{ij\alpha\beta}$.
Then the matrix $ (F_{ij\alpha\beta})$ is positive-definite
on ${\bf C^{Nr}} = {\bf C^r} \otimes {\bf C^N} $; that is,
there is a positive constant $c'$ so that

$$  \sum F_{ij\alpha\beta} t_i {\overline t_j} s_\alpha {\overline s_\beta}
\ge c' \sum |t_i s_\alpha|^2 = c' \sum |t_i|^2 \sum |s_\alpha|^2  \eqno (12) $$

Proof. First we show that 1 implies 2. Given $h \in V_d \subset {\cal A}$, we write
$h=(h_1,...,h_r)$. We assume that (6) holds and    
compute the left side of (10), to obtain

$$ \langle Th,h \rangle_{\cal H} = \int \int \sum_{k=1}^s \sum_{i,j=1}^r 
{\overline {E_{ik}(w)} } E_{jk}(z) h_i(w) {\overline {h_j(z)}} dV(w) dV(z) $$

$$ = \sum_{k=1}^s | \int \sum_j E_{jk}(z) {\overline {h_j(z)}} dV(z) |^2 = 
\sum_{k=1}^s | \langle A_k , h \rangle_{\cal H}|^2 \eqno (13) $$

The condition of strict holomorphic factorability
guarantees that the vectors $A_k$, the column vectors of $E$, form a basis for $V_d$. 
The last
expression in (13) is therefore $\ge c ||h||^2_{\cal H} $. 
This proves that 1 implies 2.

Next we show that 3 implies 1. Recall that
$N$ is the dimension of the space of homogeneous polynomials
of degree $d$ in $n$ variables, and that $(F_{ij})$ is an $r$ by $r$ matrix.
If (12) holds, there is a basis $\{ E_{i\beta}\}$
of ${\bf C^{Nr}}$ so that
$ F_{ij\alpha\beta} = \langle E_{j\alpha}, E_{i\beta} \rangle $.
Plug this in (11) to obtain

$$ F_{ij}(z,{\overline w}) = \sum F_{ij\alpha\beta} z^\alpha {\overline w}^\beta =
\sum \langle E_{j\alpha},E_{i\beta} \rangle z^\alpha {\overline w}^\beta \eqno (14) $$
Now define $A_j (z)$ by $A_j(z) = \sum E_{j\alpha} z^\alpha $. We see that

$$ F_{ij}(z,{\overline w}) = \langle A_j(z), A_i(w) \rangle \eqno (15) $$
and hence that 1 holds. 

It remains to prove that 2 implies 3. 
We write $h_i(z) = \sum H_{i\alpha} z^\alpha $, and we plug this into
$\langle Th,h \rangle_{\cal H}$. 
Recall that distinct monomials are orthogonal,
so we may write

$$ \langle Th,h \rangle_{\cal H} = 
\sum \int \int F_{ij\alpha\beta}z^\alpha {\overline w^\beta}
H_{i\mu}w^\mu {\overline H_{j\nu}} {\overline z^\nu} dV(w) dV(z) =
\sum F_{ij\alpha\beta}
H_{i\beta}{\overline H_{j\alpha}} p_\alpha p_\beta \eqno (16) $$
where the positive numbers $p_\alpha$ are equal to $||z^\alpha||^2_{L^2} $.

On the other hand, we have $||h||^2_{\cal H} = \sum |H_{i\alpha}|^2 p_\alpha$ 
by a similar
calculation. Thus (10) implies that there is a positive constant $c$ so that

$$ \sum F_{ij\alpha\beta}p_\alpha p_\beta H_{i\beta}{\overline H_{j\alpha}} 
\ge c \sum |H_{i\alpha}|^2 p_\alpha \eqno (17) $$

Since the $p_\beta$ are positive numbers, (17) implies that the
matrix with entries $ F_{ij\alpha\beta}$ is also positive-definite, with
a different constant $c'$.  This gives 3.
$\spadesuit$

\bigskip

{\bf III. Proof of Theorem 1}
\bigskip

Suppose that the entries of $F$ are bihomogeneous polynomials
of degree $2m$. Let $Q_d$ be the operator on $V_{m+d}$ whose kernel is given by 
$\langle z,w \rangle^d F(z,{\overline w})$.
In order to prove Theorem 1, Proposition 1
implies that we must find an integer $d$ so that
$Q_d$ is positive-definite on $V_{m+d}$.
We observe immediately, that 
if this holds for some $d$, then it holds for all larger
integers. See [CD1], whose title suggests this stabilization process.
Furthermore, the operators $Q_d$ are zero except on $V_{m+d}$.
This suggests considering their sum, weighted by positive
constants, on the whole space ${\cal A}$.
 
If we choose the positive constants $C_d$ appropriately, then
$$ \sum C_d \langle z,w \rangle^d = 
{n! \over \pi^n} {1 \over (1 - \langle z,w \rangle)^{n+1}} = B(z,w)  \eqno(18) $$
Here $B(z,w)$ is the Bergman kernel function for
the unit ball $B$ in complex Euclidean space ${\bf C^n}$.
The crucial property of the Bergman kernel function is that it is the integral
kernel of the Bergman projection mapping $L^2(B)$ to its closed subspace 
$A^2(B)$ of holomorphic
functions. The kernel function satisfies
$$ B(z,w) = \sum_\alpha \phi_\alpha (z) {\overline  {\phi_\alpha (w)}} $$
where the collection $\{ \phi_\alpha \}$ is any complete orthonormal set 
for the Hilbert space $A^2(B)$. 
For the unit ball, one can choose $\phi_\alpha = c_\alpha z^\alpha$, where
$c_\alpha$ is a normalizing constant, and $z^\alpha$ denotes the indicated
monomial.

Recall that ${\cal H}$ denotes the Cartesian product
of $r$ copies of $L^2(B)$. Let $P:{\cal H} \to {\cal A}$
denote the Bergman projection, acting componentwise.
Motivated by (18),
we let $Q$ denote the integral operator on ${\cal H}$ whose kernel is given by
$B(z,w) F(z,{\overline w}) $. For a scalar function $\psi$,
we let $M_\psi$ denote the operator on ${\cal H}$ given 
by multiplication by $\psi$.
Also $M_F$ denotes
matrix multiplication by $F$. There is no integration involved in these operators.
Choose a smooth non-negative
function $\phi$ of compact support that is 
positive near the origin. 

We may write

$$ Q = (M_F P + P M_\phi) + (Q - M_F P) - PM_\phi  = T_1 + T_2 + T_3 \eqno (19) $$
We claim that $T_1$ is positive, and that $T_2$ and $T_3$ are
compact. This will show that $Q = S + K$, where $S$ is
positive on ${\cal A}$ and $K$ is compact.

{\bf Lemma 1}. $T_1$ is positive on all of ${\cal A}$,
and $T_3$ is compact on ${\cal H}$.

Proof. The second statement is immediate, because the integral kernel is
smooth everywhere on the ball. The first statement follows because $P$ is
a self-adjoint projection. To see this, let $h \in {\cal A}$.

$$ \langle T_1 h, h \rangle_{\cal H} = \langle M_F P h + P M_\phi h, h \rangle_{\cal H} 
= \langle M_F h,h \rangle_{\cal H} + \langle M_\phi h, h \rangle_{\cal H} = 
\langle M_{F+\phi} h, h \rangle_{\cal H} \eqno (20) $$

Since the multiplier $F + \phi$ is strictly positive-definite at all points, 
the last expression in (20) is
at least $C ||h||^2_{\cal H} $, and the result follows. $\spadesuit$

{\bf Lemma 2}. $T_2$ is compact.

Proof. This follows from Theorem 1 in [CD2], but is elementary in this
case, because of the explicit nature of the Bergman kernel. The kernel of
$Q - M_F P$ is 
$$ {n! \over \pi^n} 
{(F(z,{\overline w}) - F(z,{\overline z})) 
\over (1 -\langle z,w \rangle)^{n+1}}  $$
The numerator (a matrix of polynomials)
vanishes on the boundary diagonal, where the only singularities of the
denominator occur. One can use Young's inequality to verify that $T_2$ is
compact. $\spadesuit$

We summarize what we have proved so far.
The operator $Q$ on ${\cal H}$ has kernel given by
$B(z,w) F(z,{\overline w})$. By Lemmas 1 and 2, we have written $Q=S+K$,
where $S$ is positive on ${\cal A}$ and $K$ is compact.

The operator $Q$ vanishes off ${\cal A}$,
and we have ${\cal A} = \oplus V_j$.
Write $Q_d$ for the restriction of $Q$ to $V_{m+d}$. 
If we show that $Q_d$ is positive on $V_{m+d}$ for sufficiently large
$d$, then an application of  Proposition 1 completes the proof
of Theorem 1.

Since $S$ is positive, there is $c>0$ so that 
$\langle Sh,h\rangle_{\cal H} \ge c ||h||^2_{\cal H}$. 
Since $K$ is compact, there is a finite rank operator $L$ such that the operator norm
$ |||K-L||| < {c \over 3} $. See [R].
Write $Q = S + L + (K-L) $ so that 

$$ \langle Qh,h \rangle_{\cal H} = \langle Sh,h \rangle_{\cal H}
+\langle Lh,h \rangle_{\cal H} + \langle (K-L)h,h \rangle_{\cal H}  \eqno (21)$$
Using the lower bound on $S$, and because 
$ |\langle (K-L) h,h\rangle_{\cal H}| \le {c \over 3} ||h||^2_{\cal H}$, we can
write 

$$ \langle Qh,h \rangle_{\cal H} \ge c||h||^2_{\cal H} - {c \over 3} ||h||^2_{\cal H}
- |\langle Lh,h\rangle_{\cal H}| \ge {2c \over 3} ||h||^2_{\cal H} -
|\langle Lh,h\rangle_{\cal H}|  \eqno (22) $$

Because $L$ is finite rank, we can choose
$d_0$ sufficiently large such that, for $d \ge d_0$,
 the restriction of $L$ to $V_{m+d}$ satisfies
$ |\langle Lh,h\rangle_{\cal H}| \le {c \over 3} ||h||^2_{\cal H}$ also.
Combining this with (22) implies that
the restriction of $Q$ to $V_{m+d}$ is positive. 
By Proposition 1, its kernel 
$ \langle z,w \rangle^d F(z,{\overline w})$ can be written
$A(w)^* A(z)$, completing the proof of Theorem 1.
$\spadesuit$

\bigskip

{\bf IV. Examples and Applications}.

\bigskip

The integer $d$ in Theorem 1
can be arbitarily large even when $F$
has fixed degree. The example
$f_c (z,{\overline z}) = |z_1|^4 + c |z_1 z_2|^2 + |z_2|^4 $ 
is positive away
from the origin for $c>-2$. 
By Theorem 1, for each $c$ with $c>-2$,
there is a minimum $d_c$ for which (3) holds. It is elementary
to show that
$d_c \to \infty$ as $c \to -2$. See [CD1].

Because the integer $d$ can be arbitarily large, the 
holomorphic mapping $A$ from
Theorem 1 can have an arbitrarily large number of components.
This fact has consequences for 
proper holomorphic mappings between balls in
different dimensions. For example, in [CD1], Theorem 1 is used to prove the
following. Given a holomorphic
polynomial $q:{\bf C^n} \to {\bf C}$
 that doesn't vanish on the closed unit ball,
there is an integer $N$ and
a holomorphic polynomial mapping $p:{\bf C^n} \to {\bf C^N}$ 
such that
${p \over q}$ is reduced to lowest terms, 
and defines a proper map between balls.
The integer $N$ can be arbitarily large. This is in sharp contrast to the
case when $n=1$, 
where the result is trivial, and we can take $N=1$ as well.

Next we give an application to symbols of differential operators.
Let $D$ be a linear partial differential operator
on real Euclidean space ${\bf R^{2n}}$ of even order $2m$.
Recall (See [F] for example) that the principal symbol, or characteristic form, $p(\xi)$
of $D$ governs whether it is elliptic. We suppose that the principal symbol
has constant coefficients.
Thus $p(\xi) = \sum_{|\alpha|=2m} c_\alpha \xi^\alpha$, and the operator is elliptic
precisely when $p$ vanishes only at the origin.  If we make the usual
identification of ${\bf R^{2n}}$ with ${\bf C^{n}}$, then we can express $D$ in terms of
the operators $ {\partial \over \partial z_j}$ and  
$ {\partial \over \partial {\overline z}_j}$. Using multi-index notation we can then
write the principal symbol as

 $$ \sum_{|\alpha| + |\beta| = 2m} c_{\alpha \beta} 
({\partial \over \partial z})^\alpha ({\partial \over \partial {\overline z}})^\beta
= f({\partial \over \partial z},{\partial \over \partial{\overline z}}). $$

In general, $f$ is not bihomogeneous. 
There are simple simple necessary and sufficient condition for a real-valued
polynomial $f(z, {\overline z})$ to be bihomogeneous of degree $2m$. 
One is that it be both homogeneous of degree $2m$ over ${\bf R}$
and invariant under replacing $z$ by $e^{i\theta}z$. (Here $e^{i\theta}$ is a
scalar, not an $n$-tuple). Another is the existence of 
holomorphic polynomial mappings $A$ and $B$, each
homogeneous of degree $m$, such that
$ f(z, {\overline z}) = ||A(z)||^2 - ||B(z)||^2 $.
An arbitrary real-valued polynomial $p$ can be written 
as the difference of squared norms
of holomorphic polynomials, but the holomorphic polynomials will not 
be homogeneous of the same degree when $p$ fails to be bihomogeneous.
See [D] for uses of such a decomposition.

We say that the partial differential operator 
$D$ on ${\bf R^{2n}}$ is
{\it complex bihomogeneous} if its principal symbol is a bihomogeneous
polynomial. In this case we may apply Theorem 1
to obtain the following conclusion. We write $\triangle$ for the Laplace operator
defined by 
$\sum {\partial \over \partial z_j} {\partial \over \partial {\overline z}_j} $.

{\bf Theorem 2}. Let $D$ be a complex bihomogeneous linear partial differential
operator. Suppose that $p$ is the absolute value of the principal symbol of $D$.
Let $q_d$ be the absolute value of the principal symbol of 
$\triangle^d D$.
The following are equivalent:

1) $D$ is elliptic (that is, $p(z,{\overline z}) > 0$ for $z \ne 0$).

2) There is an integer $d$ and a positive-definite matrix $(E_{\mu \nu})$
so that $q_d$ satisfies 

$$  q_d ({\partial \over \partial z} ,{\partial \over \partial {\overline z}}) = 
\sum_{|\mu|=m} \sum_{ |\nu| = m} E_{\mu \nu} 
({\partial \over \partial z})^\mu ({\partial \over \partial {\overline z}})^\nu. $$

3) There is an integer $d$ so that $q_d$ is a squared norm
of a holomorphic differential operator:

$$ q_d ({\partial \over \partial z} ,{\partial \over \partial {\overline z}}) =
|| \sum A_\mu ({\partial \over \partial z})^\mu||^2 
= \sum_i |\sum A_{\mu i} ({\partial \over \partial z})^\mu |^2.
\eqno (23)$$
We assume also that the indicated
homogeneous polynomials span $V_{m+d}$.

4. There is an integer $d'$ so that $q_{d'}$ 
satisfies (23), and such that
${\bf V}(A) = \{ 0 \}$.

{\bf Proof}. The principal symbol of $\triangle^d D$, evaluated 
at $(z, {\overline z})$,  is
$q_d(z,{\overline z}) =||z||^{2d} p(z, {\overline z})$. The operator $D$ is elliptic
precisely when $|p|$ is strictly positive away from the origin.
Therefore, by Theorem 1,
$D$ is elliptic if and only there is
$d$ so that $||z||^{2d} p(z, {\overline z})$ satisfies
any of the equivalent conditions of Proposition 1.
Equation (10), applied when $r=1$, is equivalent to
the positive-definiteness of
$(E_{\mu \nu})$. The strict holomorphic factorability there
is equivalent to statement 3 here.
Statement 3 obviously implies statement 4, which in turn implies
that $p(z,{\overline z})$ is positive away from the origin.
Thus the four statements are equivalent. $\spadesuit$

{\bf Remark}. This result extends to systems of
PDE in a straightforward fashion.

\bigskip

{\bf V. Reinterpretation of Theorem 1}
\bigskip

Next we 
reinterpret Theorem 1 in terms of pullbacks of
the universal
bundle over Grassman manifolds. See [W] for more details
about the universal bundle.
Let ${\bf G}_{p,N}$ denote the Grassman manifold
of $p$ planes in complex $N$-space. 
When $p=1$ we have complex projective space, 
and we write as usual ${\bf P}^{N-1}$ for 
${\bf G}_{1,N}$. Let ${ \bf U}_{p,N}$ 
denote the universal bundle over 
 ${\bf G}_{p,N}$. This bundle is sometimes known as the
tautological bundle or the stupid bundle;
a point in  ${\bf U}_{p,N}$ is a pair
$(S,\zeta)$ where $S$ is a $p$-dimensional subspace of
${\bf C^N}$ and
$\zeta \in S $. 

We let $g_0$ denote the Euclidean metric on
${\bf U}_{p,N}$. In terms of a local
frame $e$ of $ {\bf U}_{p,N}$, we define
$ g_0 (e) = e^* e $. 
Observe that if $T$ determines a change of
frames by acting on the right, then 
$$ g_0(eT) = (eT)^* (eT) = T^* e^* e T =  T^* g_0(e) T .$$
This is the correct transformation law, and hence $g_0$
defines a Hermitian metric on  $ {\bf U}_{p,N}$.
We may consider the matrix representation
of $g_0$, with respect to a local frame.
We have $(g_0)_{ij} =  \langle e_j, e_i \rangle $ 
where $\langle,\rangle$
denotes the usual Hermitian inner product on complex Euclidean
space ${\bf C^N}$, and the vectors $e_i$ 
for $i=1,...,r$ are linearly independent.
Note the interchange of indices.
We see immediately
that $g_0$ is of the form $A^* A$.

Let $L$ denote
the universal line bundle $L = {\bf U} = {\bf U}_{1,n}$
over complex projective space
${\bf P}^{n-1}$. We consider also
its $m$-th tensor power ${\bf U}^m$. Let $E$ denote
the vector bundle over ${\bf P}^{n-1}$
equal to the direct sum of $r$ copies of
${\bf U}^m$. On $L$ we use the Euclidean metric, written
$||z||^2$, and on $E$ we use the metric determined by a
matrix of bihomogeneous polynomials $F(z,{\overline z})$ that is
positive-definite for $z \ne 0$. Theorem 1 now admits the
following restatement.

{\bf Theorem 3}. Suppose that $L$ and $E$ are the bundles over
${\bf P}^{n-1}$ as described in the previous paragraph, equipped
with the given metrics. Then there are integers $N$ and
$d$, so that
the bundle $ L^d \otimes E$ over ${\bf P}^{n-1}$,
with metric determined by 
$||z||^{2d} F(z,{\overline z})$, is the isometric pullback
of the vector bundle
$ {\bf U}_{r,N}$, with the Euclidean metric,
 over the Grassmanian ${\bf G}_{r,N}$ via a holomorphic
embedding.

The link to bihomogeneous polynomials arises because
one can identify homogeneous polynomials of degree
$m$ on ${\bf C^n}$ with
sections of the $m$-th
power $H^m$ of the hyperplane bundle $H$ over ${\bf P}^{n-1}$. 
The bundle $H$ is dual to ${\bf U}$. A
matrix of bihomogeneous
polynomials determines a metric on the direct sum of $r$ copies
of ${\bf U}$.

The general result in [CD3] considers certain
base manifolds $M$, a line bundle $L$ and a vector bundle $E$
over $M$, and metrics $R$ and $F$
on them satisfying certain conditions. The
conclusion again guarantees the existence of integers $N$ and
$d$ so that $ L^d \otimes E$, with metric $R^d F$, 
is the isometric pullback
of ${\bf U}_{r,N}$, with the Euclidean metric,
 over the Grassmanian ${\bf G}_{r,N}$. The proof again relies on
the Bergman kernel and facts about compact operators, but it is
technically more difficult than the special case considered
here. 

\bigskip

{\bf VI. Remarks on classical factorization}
\bigskip

We briefly mention some of the results in [Dj] and [RR]. 
Djokovic [Dj] considers for example an $r$ by $r$
positive semi-definite matrix $F(\lambda, \mu)$ whose 
entries are complex-valued homogeneous polynomials
of degree $2m$ in the pair of real variables $(\lambda,\mu)$.
He proves that one can write
$F(\lambda,\mu) = A(\lambda,\mu)^* A(\lambda,\mu)$ where the entries in $A$ are
homogeneous polynomials of degree $m$. Two nice things about
this result are that it is not required to multiply
$F$ by powers of a scalar function, 
and $F$ is allowed to be semi-definite. On the other
hand, the theorem holds only when
the entries depend upon two real variables, the analogue of one complex
variable. For us, making $A(z)$ depend holomorphically on $z$ in (4)
requires that we work with bihomogneous polynomials. 
The only bihomogeneous polynomials
in one complex variable are constants times $|z|^{2m}$. 
Hence we could factor this scalar
out of the matrix completely, and the 
analogue of the result in [Dj] becomes trivial in our setting.
The idea of multiplying by powers of a scalar factor does not appear in [Dj].

The work in [RR] concerns functions from either ${\bf R}$ or the unit circle $S^1$ 
that
take values in non-negative operators on a Hilbert space. 
The authors study many aspects of the
factorization question in detail, including holomorphic
extension to the upper half plane or to the
unit disc.  One of many results there is that if $P(x)$ is a non-negative
operator on a Hilbert space, that is a polynomial of degree $2m$ in the real variable
$x$,
then there is an operator $Q(x)$ such that $P(x) = Q(x)^* Q(x) $ and such that $Q$ is a
polynomial of degree $m$ in $x$. Other results in [RR] are related to an
application of Theorem 1 here from [CD1]. Suppose that $f(z,{\overline z})$ is an
arbitrary polynomial that is positive on the unit
sphere. Then there is a holomorphic
polynomial mapping $g$ such that $f(z,{\overline z}) = ||g(z)||^2$ on the unit sphere.
In [RR] however positivity questions are considered only for functions depending on
one variable.

\bigskip

{\bf References}

\bigskip

[CD1] David W. Catlin and John P. D'Angelo,
A stabilization theorem for Hermitian forms and applications to holomorphic
mappings, Math Research Letters 3 (1996), 149-166. 

[CD2] David W. Catlin and John P. D'Angelo, Positivity conditions for bihomogeneous
polynomials,
Math Research Letters 4 (1997), 1-13.

[CD3] David W. Catlin and John P. D'Angelo,
Isometric embeddings of bundles, (in preparation).

[D] John P. D'Angelo, Several complex variables and the geometry of 
real hypersurfaces,
CRC Press, Boca Raton, 1993.

[Dj] D. Z. Djokovic, Hermitian matrices over polynomial rings, J.
Algebra 43 (1976), 359-374.

[F] Gerald B. Folland, Introduction to Partial Differential Equations,
Princeton University 
Press, 1976.

[RR] M. Rosenblum and J. Rovnyak, The factorization problem for
nonnegative operator valued functions, Bulletin A.M.S. 77 (1971),
287-318.

[Ru] Walter Rudin, Functional Analysis, McGraw-Hill, New York, 1973.

[W] Raymond O. Wells, Differential Analysis on Complex Manifolds, 
Prentice-Hall,
Englewood Cliffs, New Jersey, 1973.

\end